\documentclass[conf]{new-aiaa}
\usepackage[utf8]{inputenc}

\usepackage[english]{babel}
\usepackage{blindtext}

\usepackage{graphicx}
\usepackage[version=4]{mhchem}
\usepackage{siunitx}
\usepackage{longtable,tabularx}
\usepackage{color,url}

\newcommand{\mbf}{\mathbf}

\title{A Koopman-Operator Control Optimization for Relative Motion in Space}

\author{Simone Servadio\footnote{Postdoctoral Associate, Department of Aeronautics and Astronautics, Massachusetts Institute of Technology, MA 02139, USA. email: simoserv@mit.edu.}}
\affil{Massachusetts Institute of Technology, Cambridge, MA, 02139, USA}
\author{Roberto Armellin\footnote{Professor, Mechanical Engineering, roberto.armellin@auckland.ac.nz}}
\affil{Auckland Space Institute, The University of Auckland, Auckland, NZL}
\author{ Richard Linares\footnote{Boeing Assistant Professor, Department of Aeronautics and Astronautics, Massachusetts Institute of Technology, MA 02139, USA. email: linaresr@mit.edu}}
\affil{Massachusetts Institute of Technology, Cambridge, MA, 02139, USA}

\begin{document}

\maketitle

\begin{abstract}
A high order optimal control strategy implemented in the Koopman operator framework is proposed in this work. The new technique exploits the Koopman representation of the solution of the equations of motion to develop an energy optimal inverse control methodology. The operator theory can reformulate a nonlinear dynamical system of finite dimension into a linear system with an infinite number of dimensions. As a results, the state of any nonlinear dynamics is represented as a linear combination of high-order orthogonal polynomials, which creates the state transition polynomial map of the solution. Since the optimal control technique can be reduced to a two-points boundary value problem, the Koopman map is used to connect the state and control variables in time, such that optimal values are obtained through map inversion and polynomial evaluation. The new technique is applied to rendezvous applications in space, where the relative motion between two satellites is modelled with a high-order polynomial series expansion of the Lagrangian of the system, such that the Clohessy-Wiltshire equations represent the reduction of the high-order model to a linear truncation. 
\end{abstract}

\section{Introduction}
Trajectory optimization is a key element in every space mission \cite{conway2010spacecraft}. Trajectory optimization enables complex interplanetary transfers \cite{jehn2012trajectory}, operative orbit acquisition (e.g., from a geostationary transfer orbit to a geostationary one \cite{haberkorn2004low}), rendezvous and docking \cite{luo2014survey}, formation flying control and reconfiguration \cite{di2018survey} and station-keeping \cite{gazzino2019long}, just to mention a few. The minimization of propellant or transfer time is the typical objective of trajectory optimization.  

Low-thrust propulsion is more frequently used as it allows to reduce the propellant mass thanks to high specific impulse. Trajectory optimization with low-thrust propulsion is a challenging task as it involves solving complex optimal control problems (OCP) \cite{betts1998survey}.  Two main methods are used to solve low-thrust optimization problems: direct and indirect. In the direct methods, the control is parametrized and the dynamics are reduced into a set of constraints. In this way, the OPC becomes a parametric optimization problem solved with nonlinear programming tools. The large dimension of the problem and possible convergence to suboptimal solutions are two main issue of direct methods.
Indirect methods can find optimum solutions more accurately and result in a lower dimension numerical problem. These methods apply calculus of variation and Pontryagin’s maximum principle to transform the OCP into a two-point boundary value problem \cite{bryson2018applied}. The main drawback of indirect methods is that the convergence
domain is small, guesses for initial costate are not easily found, and discontinuity in the control laws are difficult to handle \cite{jiang2012practical}.

In the astrodynamics community there is a strong interest in developing techniques to solve OCP in an efficient way for set of initial and final conditions. This interest is fueled by the need for performing trade-off analyses efficiently, quantifying the effect of uncertainties and better estimate mission budgets, or use OCP solutions for spacecraft guidance. The formulation direct methods as convex optimization problems has received a lot of attention in recent years \cite{malyuta2021advances}. Thanks to the convergence properties of convex optimization and the efficiency of modern solvers, OCP can be potentially solved onboard and the solutions used for spacecraft guidance \cite{armellin2021collision,hofmann2021rapid}.  

An alternative approach for autonomous guidance is to computing the expansion of the solution of OCPs with respect to uncertainties in initial and final conditions, using state transition tensors (STTs) \cite{boone2021orbital,boone2021variable} and differential algebra (DA) \cite{di2014high,di2014high2}. These expansions can allow for the update of the controls due to the presence of errors in the spacecraft or target states. Indeed the required update of the reference trajectory is reduced to a simple evaluation of polynomials.

In the same vein, this work introduces an alternative way to efficiently compute the solution of OCP formulated in a indirect approach. Our approach is based on the Koopman operator (KO). The KO theory was initially developed by Koopman \cite{koopman1931hamiltonian} and further developed by Von Neumann in the 1930s \cite{neumann1932operatorenmethode}. They transitioned the operator theory initially developed in quantum mechanics \cite{naylor2000,prigogine} to classical mechanics problems. The key result was the development of a linear operator that can reformulate a nonlinear dynamical system of finite dimension into a linear system with an infinite number of dimensions. The linearity of the KO is very appealing, but this benefit is contrasted with the fact that it is infinite-dimensional. However, this issue can be overcome by capturing the evolution on a finite subspace spanned by a finite set of basis functions. In effect, this is a truncation of the KO to a finite subspace. More recently, the KO has been employed in astrodynamics to model the motion of a satellite about an oblate planet~\cite{arnas2021approximate,arnas2021analysis} and to analyze the orbits of a satellite in the close proximity of libration points~\cite{servadio2021dynamics}.

The KO solution of the system represents the time propagation of the state as a linear combination of well selected polynomial basis functions. This works exploit this representation to create a high-order state polynomial transition map to solve the TPBVP, following the inverse control approach. Indeed, the map can be inverted~\cite{di2014high} to obtain the optimal control that connects, in the desired amount of time, the state of the system between the boundary conditions. The final result is a high-order inverse control technique that does not require the evaluation of a reference solution, contrary to previous approaches \cite{di2014high,di2014high2}.

The problem we aim to solve with the KO-based approach is the one of a rendezvous with a chief spacecraft moving in a circular Keplerian orbit. The transfer is implemented with a low-thrust propulsion system and the objective is to minimize the energy of the transfer. Differently from the classical Clohessy-Wiltshire equations \cite{curtis2013orbital}, we retain terms up to an arbitrary high order in the Legendre expansion of the  gravitational potential. The equations of motion are derived from the polynomial series expansion of the Lagrangian of the system \cite{kasdin2005canonical}, where the accuracy of the representation of relative motion increases due to the presence of high order terms \cite{curtis2013orbital}. This results in a nonlinear control problem that can only be solved numerically due to the nonlinearities in the Euler-Lagrange equations. 

The paper is structured as follows. The Koopman operator theory is presented in the following section. Afterwards, the energy optimal control technique is derived in the KO framework, highlighting the advantages of the newly developed methodology. Then, the numerical application are presented. First, the Duffing oscillator is controlled to give a representation of the new technique for a easy problem. Secondly, the linear and high-order controlled Clohessy-Wiltshire equations are derived, providing comparison between the two dynamics. The optimal inverse KO control is applied to achieve rendezvous between a chaser spacecraft and a target for different times of flight and boundary conditions. Lastly, conclusions are drawn. 

%%%%%%%%%%%%%%%%%%%%%%%%%%%%%%%%%%%%%%%%%%%%%%%%%%%%%%
%%%%%%%%%%%%%%%%%%%%%%%%%%%%%%%%%%%%%%%%%%%%%%%%%%%%%%
%%%%%%%%%%%%%%%%%%%%%%%%%%%%%%%%%%%%%%%%%%%%%%%%%%%%%%
%%%%%%%%%%%%%%%%%%%%%%%%%%%%%%%%%%%%%%%%%%%%%%%%%%%%%%

\section{The Koopman Operator} \label{KOsec}
A classical definition of nonlinear dynamical systems is given by the initial value problem, which can be represented by a set of coupled autonomous ordinary differential equations in the form:
\begin{equation}\label{problem}
\left\{ \begin{tabular}{l}
    $\displaystyle\frac{d}{dt}{\bf x}(t) = {\bf f}({\bf x})$ \\
    ${\bf x}(t_0) = {\bf x_0}$ 
\end{tabular}  \right.
\end{equation}
where ${\mbf x}\in\mathbb{R}^m$ is the state  which depends on the time evolution $t$, ${\bf f}:\mathbb{R}^m\rightarrow \mathbb{R}^m$ is the nonlinear dynamics model, $m$ is the number of dimensions in which the problem is defined, and ${\bf x}_0$ is the initial condition of the system at time $t_0$. The Koopman operator (KO), $\mathcal{K}(\cdot)$ is an infinite-dimensional linear operator that evolves all observable functions $\mathcal{G}({\bf x})$ of the state and it allows to redefine any problem written in classical mechanics into an equivalent Koopman formulation.

Let $\mathcal{F}$ be a vector space of observable functions, where $\mathcal{G}({\bf x})\in \mathcal{F}$. Since the KO is an infinite-dimensional linear operator, this space of functions $\mathcal{F}$, which the observables are defined on, is also infinite-dimensional. Therefore, if $g\subseteq\mathcal{G}({\bf x})$ is a given observable in this space, the evolution of $g$ in the dynamical system is represented by:
\begin{equation}
\mathcal{K}\left(g({\bf x})\right) = \frac{d}{dt}g({\bf x}) = \left( \nabla_{{\bf x}} g({\bf x})\right)\frac{d}{dt}{\bf x}(t) = \left( \nabla_{{\bf x}} g({\bf x})\right){\bf f}({\bf x}),
\end{equation}
where $\nabla_{{\bf x}} g = (\partial g/\partial x_1,\partial g/\partial x_2,\dots,\partial g/\partial x_m)$. That way, the evolution of any observable subjected to the dynamical system is provided by the Koopman operator:
\begin{equation}
\mathcal{K}\left(\cdot\right) = \left( \nabla_{{\bf x}} \cdot\right){\bf f}({\bf x}).
\end{equation}
Note that the evolution of the observables is provided by the application of the chain rule for the time derivative of $g({\bf x})$. Consequently, the defined operator is linear~\cite{koopman1931hamiltonian}, such that:
\begin{equation}
    \mathcal{K}\left(\beta_1g_1({\bf x})+\beta_2g_2({\bf x})\right)=\beta_1\mathcal{K}\left(g_1({\bf x})\right)+\beta_2\mathcal{K}\left(g_2({\bf x})\right),
\end{equation}
for any pair of observables $g_1\subseteq\mathcal{G}({\bf x})$ and $g_2\subseteq\mathcal{G}({\bf x})$ and any arbitrary constants $\beta_1$ and $\beta_2$. The linearity of the Koopman Operator is very appealing, but this benefit is contrasted with the fact that it is infinite-dimensional. However, this issue can be overcome by capturing the evolution of the system on a finite subspace spanned by a finite set of basis functions instead of all measurement functions in a Hilbert space. In effect, this is a truncation of the Koopman operator to a finite subspace $\mathcal{F}_D$ of dimension $n$, where $\mathcal{F}_D \in \mathcal{F}$. This subspace $\mathcal{F}_D$ can be spanned by any set of eigenfunctions $\phi_i\in\mathcal{F}_D$, with $i\in\{1,2,\dots,n\}$, defined as:
\begin{equation}\label{koopman}
\mathcal{K}\left(\phi_i({\bf x})\right) =\frac{d}{dt}\phi_i({\bf x})=\lambda_i \phi_i({\bf x}),
\end{equation}
where $\lambda_i$ are the eigenvalues associated with the eigenfunctions $\phi_i$, and $n$ is the number of eigenfunctions chosen to represent the space. Therefore, the Koopman eigenfunctions can be used to form a transformation of variables that linearizes the system. Particularly, let ${\bf \Phi}({\bf x}) = \left(\phi_1({\bf x}), \dots, \phi_n({\bf x}) \right)^T$ be the set of eigenfunctions of the KO in $\mathcal{F}_D$. Then, using the relation in Eq.~\eqref{koopman}, it is possible to write the evolution of ${\bf \Phi}$ as: 
\begin{equation}\label{time_eigenf}
\mathcal{K}\left(\bf \Phi\right) = \frac{d}{dt}{\bf \Phi}=\Lambda {\bf \Phi},
\end{equation}
where $\boldsymbol \Lambda=\text{diag}([\lambda_1, \dots,\lambda_n])$ is the diagonal matrix containing the eigenvalues of the system in $\mathcal{F}_D$. This transformation is called the Koopman Canonical Transform~\cite{surana2016linear}. The solution of Eq.~\eqref{time_eigenf} is:
\begin{equation}\label{eigen_time}
    {\bf \Phi}(t) = \exp(\boldsymbol\Lambda t){\bf \Phi}(t_0),
\end{equation}
where ${\bf \Phi}(t_0)$ is the value of the eigenfunctions at the initial time $t_0$. This result will be used later to solve the complete system once the eigenfunctions of the operator are obtained.

In general, we are interested in the identity observable, that is, ${\bf g}({\bf x})={\bf x}$. Therefore, it is required to be able to represent these observables in terms of the KO eigenfunctions. This is achieved using the Koopman modes, i. e., the projection of the full-state observable onto the KO eigenfunctions. If this projection can be found, the evolution of the state is represented by means of the evolution of the KO eigenfunctions, and thus, an approximate solution to the system can be provided. 

%%%%%%%%%%%%%%%%%%%%%%%%%%%%%%%%%%%%%%%%%%%%%%%
%%%%%%%%%%%%%%%%%%%%%%%%%%%%%%%%%%%%%%%%%%%%%%%

\subsection{The Galerkin Method}
The eigenfunction of the Koopman operator are obtained through the Galerkin method. The eigenfunctions give rise to a set of linear first-order PDEs for the eigenfunctions in the form: 
\begin{equation}\label{koopman2}
\mathcal{K}(\phi_i) = \left(\nabla_x \phi_i({\bf x})\right) {\bf f}({\bf x}) =\lambda_i \phi_i({\bf x}),
\end{equation}
or in a more expanded notation:
\begin{equation}\label{koopman_PDE}
\frac{d \phi_i({\bf x})}{dt}={ f}_1\left({\bf x}\right)\frac{\partial }{\partial x_1} \phi_i({\bf x}) +\cdots+  { f}_d\left({\bf x}\right)\frac{\partial }{\partial x_d}\phi_i({\bf x}) =\lambda_i \phi_i({\bf x}),
\end{equation}
where ${\bf f}({\bf x}) = (f_1({\bf x}), f_2({\bf x}), \dots, f_m({\bf x}))^T$. This equation is a linear first-order PDE and in general has no closed-form solution. However, it is possible to approximate the solution using the Galerkin method.

This work makes use of Legendre polynomials, due to the advantages they provide in the computation of the Koopman matrix \cite{arnas2021approximate}. Legendre polynomials are a set of orthogonal polynomials defined in a Hilbert space that generate a complete basis. The idea of this methodology is to represent any function of the space by using this set of basis functions. This is done by the use of inner products and the correct normalization of the Legendre polynomials. Let $f$ and $g$ be two arbitrary functions from the Hilbert space considered. Then, the inner product between these two functions is defined as:
\begin{equation}
\langle  f, g \rangle =\int_{\Omega} f({\bf x})g({\bf x}) w({\bf x})d{\bf x} \label{innerproduct},
\end{equation}
where $w({\bf x})$ is a positive weighting function defined on the space domain $\Omega$. This work makes use of the normalized Legendre polynomials as proposed by Arnas and Linares~\cite{zonal_koopman}, due to the advantages they provide in the computation of the inner product and in the accuracy of obtained in the Koopman matrix. First, normalized Legendre polynomials are a set of orthonormal polynomials, which simplifies the projection process. Second, Legendre polynomials are defined in a bounded domain $\Omega = [-1,1]$, as opposed to other basis functions, like Hermite polynomials, that are defined in all $\mathbb R$. Third, their weighted function is a constant, $w({\bf x}) = 1$. This aspect not only simplifies the expression of the inner product, such that it becomes easier to implement a convenient automatic computation of the integrals, but also approximates functions uniformly in the domain of definition~\cite{zonal_koopman,servadio2021koopman}. In addition, the normalized Legendre polynomials are defined such that:
\begin{equation} \label{norm}
\langle \mathcal L_i, \mathcal L_j \rangle =\int_{\Omega} \mathcal L_i({\bf x}) \mathcal L_j({\bf x}) w({\bf x})d{\bf x} = \delta_{ij},
\end{equation}
where $\mathcal L_i$ and $\mathcal L_j$ with $\{i,j\}\in\{1,\dots,n\}$ are two given normalized Legendre polynomials from the set of basis functions selected, and $\delta_{ij}$ is Kronecker's delta.

%%%%%%%%%%%%%%%%%%%%%%%%%%%%%%%%%%%%%%%%%%%%%%%
%%%%%%%%%%%%%%%%%%%%%%%%%%%%%%%%%%%%%%%%%%%%%%%

\subsection{Time Evolution of the Basis Functions and the Koopman Matrix}
The KO is based on the selection of a set of $n$ orthogonal basis functions that describe the space of solutions. Legendre polynomials have been selected as a base for the theory developed in this paper. Therefore, following the KO approach, it is desired to describe the rate of change of the basis functions, in time, as a linear combination of the functions themselves:
\begin{equation}
    \dfrac{d\mathcal L}{d t} = {\mbf K} \mathcal L \label{eqlin}
\end{equation}
where $\mathcal L$ are the Legendre polynomials and ${\mbf K} $ is the Koopman matrix. Equation (\ref{eqlin}) is linear but high dimensional. The original ODE is expressed as a combination of linear differential equations, which solution would be exact for an infinite dimensional solution space. However, due to practicality, the Koopman solution of the ODE is an approximation that works on a subset of the infinite Hilbert space that correctly describes the solution. As such, the higher the order of the KO, the more accurate the Koopman solution, since the observables are described using more eigenfunctions that better deal with high nonlinearities. Each row of Eq.~(\ref{eqlin}) describes how each single total derivative of the Legendre polynomial behaves as a linear combination of the Legendre polynomials themselves. Therefore, each entry of the Koopman matrix is evaluated through the inner product, which calculates the projection of the derivative into the basis functions. Thus, the $i$th row can be written as
\begin{align}
\dfrac{d \mathcal L_i}{d t} &= \langle \dfrac{d \mathcal L_i}{d t} , \mathcal L_0\rangle \mathcal L_0 +   \langle \dfrac{d \mathcal L_i}{d t} , \mathcal L_1\rangle  \mathcal L_1 +   \langle \dfrac{d \mathcal L_i}{d t} , \mathcal L_2\rangle  \mathcal L_2 +   \dots = \sum_{j=0}^{n} \langle \dfrac{d \mathcal L_i}{d t} , \mathcal L_j\rangle \mathcal L_j \label{element}
\end{align}
From this expression, the single entry $\mbf{ K}_{ij}$ of the Koopman matrix can be evaluated  as 
\begin{equation}
\mbf{ K}_{ij}=\langle \dfrac{d \mathcal L_i}{d t} , \mathcal L_j\rangle
\end{equation}
that indicates the projection of the rate of change of the $i$th basis function onto the $j$th basis function. Therefore, the KO matrix is square with dimensions $n \times n$. This formulation completely describes the evolution of the system in time, for each basis function. The inner product has already been defined, and it evaluates each component of the Koopman matrix via the Galerkin method. The remaining term that requires handling is the total derivative of the basis functions. However, using chain rule, we can decompose each total derivative in terms of its partial derivatives with respect to the states. Considering the general case of a system with $m$ states, the $i$th total derivative is expanded as 
\begin{align}
\dfrac{d \mathcal L_i}{d t} &= \dfrac{\partial \mathcal L_i}{\partial x_1} \dfrac{d x_1}{d t} + \dfrac{\partial \mathcal L_i}{\partial x_2} \dfrac{d x_2}{d t} + \dfrac{\partial \mathcal L_i}{\partial x_3} \dfrac{d x_3}{d t} +  \dots  = \sum_{j=1}^{m} \dfrac{\partial \mathcal L_i}{\partial x_j} \dfrac{d x_j}{d t} = \sum_{j=1}^{m} \dfrac{\partial \mathcal L_i}{\partial x_j}{ f}_j({\bf x}). \label{totd}
\end{align}
By looking at the equation, it can be noted that each partial derivative is always known and it does not depend directly on the problem. On the contrary, $ { f}_j({\bf x}) $ is one equation of motion, which changes depending on the problem and that describe the selected application.

%%%%%%%%%%%%%%%%%%%%%%%%%%%%%%%%%%%%%%%%%%%%%%%
%%%%%%%%%%%%%%%%%%%%%%%%%%%%%%%%%%%%%%%%%%%%%%%

\subsection{The Koopman Solution of the System}
The system is described through the Koopman matrix. However, in order to obtain the solution of the system, we are interested in the time evolution of the state. Thanks to the Koopman representation, any function of the state can be approximated as a linear combination of the basis functions. Let $\mathbf g (\mbf x)$ be a set of observables that we are interested to evaluate. Any $i$th observable can be projected on the set of basis functions using the inner product:
\begin{equation}
    \mbf g = \sum_{j=0}^{n} \langle \mbf g, \mathcal L_j\rangle \mathcal L_j  
\end{equation}
These observables can be any function of the original variables $\mbf x$, including the states themselves, which constitutes identity observables. Therefore, the set of coefficients that expresses the projection of the observables into the basis functions can be represented in matrix form, likewise for the Koopman matrix:
\begin{equation}
    \mbf H_{ij} = \langle \mbf g_i , \mathcal L_j\rangle
\end{equation}
where $\mbf H_{i,j}$ indicates the projection of the $i$th observable onto the $j$th basis function. The number of observables dictates the dimensions of $\mbf H$: in the particular case of the identity observables, where it is desired to obtain the state of the system, matrix $\mbf H$ has dimensions $m \times n$. 

The observables are represented as a linear combination of the basis functions, which time evolution is approximated by the Koopman matrix. Under the assumption of matrix $\mbf K $ being diagonalizable, the eigendecomposition of the dynamics can be performed;  
\begin{equation}
   \mbf V \mbf K = \boldsymbol \Lambda \mbf V   \label{eig}
\end{equation}
where $\mbf V$ is the eigenvectors matrix and $\boldsymbol \Lambda$ is the diagonal matrix of eigenvalues. Rewriting Eq.~(\ref{eqlin}) by highlighting the variable dependency of the functions
\begin{equation}
    \dfrac{d}{d t} (\mathcal L(\mbf x(t)))= \mbf K \mathcal L(\mbf x(t)), \label{eqlin_dep}
\end{equation}
thanks to the eigendecomposition of the Koopman matrix, Eq.~(\ref{eig}), the eigenfunctions of the system are obtained from the Legendre polynomials using the eigenvectors matrix:
\begin{equation}
    \boldsymbol \phi(\mbf x(t)) = \mbf V \mathcal L(\mbf x(t)) \label{eigfun}
\end{equation}
The analysis of the eigenfunctions brings us to a simpler description of the system where each differential equation is decoupled from the others, and, therefore, easier to solve. Let us take the time derivative of the eigenfunctions and substitute what has been derived so far:
\begin{subequations}
\begin{align}
    \dfrac{d}{d t }\boldsymbol \phi(\mbf x(t)) &= \dfrac{d}{d t }  (\mbf V \mathcal L(\mbf x(t))) = \mbf  V \dfrac{d}{d t }  ( \mathcal L(\mbf x(t))) \\
     &= \mbf V \mbf K \mathcal L(\mbf x(t)) \\
     &= \boldsymbol  \Lambda \mbf  V \mathcal L(\mbf x(t)) \\
     &= \boldsymbol  \Lambda \boldsymbol \phi(\mbf x(t)) \label{diag}
\end{align}
\end{subequations}
where the substitutions come, respectively, from Eq.~(\ref{eqlin_dep}), Eq.~(\ref{eig}), and Eq.~(\ref{eigfun}). The last relation, Eq.~(\ref{diag}), shows a diagonal system of ODEs whose solution is known
\begin{equation}
    \boldsymbol \phi(\mbf x(t)) = \exp(\boldsymbol  \Lambda t) \boldsymbol \phi(\mbf x(t_0)) \label{sol}
\end{equation}
where $\boldsymbol \phi(\mbf x(t_0))$ indicates the value of the eigenfunctions at the initial time $t_0$. 

The evolution of the eigenfunctions with time is used to find the solution of any observable, and thus, of the state of the system (when the observable is the identity). Each observable function $\mbf g (\mbf x (t))$ has already been represented as a linear combination of the basis functions of the system through matrix $\mbf H$. Therefore, after some manipulations, the time evolution of the observables can be expressed solely as a function of time $t$ and of the state initial condition $\mbf x(t_0)$ :
\begin{subequations}
\begin{align}
   \mbf g (\mbf x (t)) &= \mbf H \mathcal L(\mbf x(t)) \\
   &= \mbf  H \mbf V^{-1} \boldsymbol \phi(\mbf x(t))\\
   &= \mbf H \mbf V^{-1} \exp(\boldsymbol\Lambda t) \boldsymbol \phi(\mbf x(t_0))\\
   &= \mbf H\mbf V^{-1} \exp(\boldsymbol\Lambda t)\mbf V \mathcal L(\mbf x(t_0))  \label{final}
\end{align}
\end{subequations}
where the substitutions come, respectively, from the inversion of Eq.~(\ref{eigfun}), Eq.~(\ref{sol}), and Eq.~(\ref{eigfun}). The solution of the equations of motion has been found, having picked the state of the system as observables. 

A more complete and detailed explanation of the Koopman framework, with a curated analysis and derivation of the mathematics of the operator, is offered in Ref.~\cite{arnas2021approximate,servadio2021koopman,arnas2021analysis,servadio2021dynamics}.

%%%%%%%%%%%%%%%%%%%%%%%%%%%%%%%%%%%%
%%%%%%%%%%%%%%%%%%%%%%%%%%%%%%%%%%%%
%%%%%%%%%%%%%%%%%%%%%%%%%%%%%%%%%%%%

\section{Optimal Control via Koopman Polynomial Map Inversion}
The KO solution of the system is expressed, for each time step, as a linear combination of the system basis functions, which get evaluated at the initial condition. This representation suits the evaluation of optimal control via indirect method. 

\subsection{Energy Optimal Costate Evaluation}
Assume that we want to link two states, $\mbf x_0$ and $\mbf x_f$ in a given amount of time, time-of-flight (ToF) $\bar t$, 
\begin{equation}
    \mbf x(t_0) = \mbf x_0, \quad \quad \mbf x(t_f) = \mbf x_f, \quad \quad t_f - t_0 = \bar t, 
\end{equation}
consuming minimum energy. The costs function $\mathcal J$ related to the performance index $\mathcal P$ to be optimized is 
\begin{equation}
    \mathcal J = \int^{t_f}_{t_0} \mathcal P dt = \int^{t_f}_{t_0} \dfrac{1}{2}\mbf \tau^2 dt
\end{equation}
After separating the state of the system into position and velocity, $\mbf x = [  \mbf r \quad   \mbf v ]^T$, the problem statement can be written as:
\begin{equation}
\mbf p(\mbf x, \tau, \hat{\boldsymbol \alpha}) = 
    \begin{cases}
      \dot{\mbf r} = \mbf v \\ 
      \dot{\mbf v} = \mbf f\big( \mbf r, \mbf v \big) + \tau \hat{\boldsymbol \alpha}  
    \end{cases} 
\end{equation}
where $\hat{\boldsymbol \alpha}$ describes the direction of the control, while $\mbf f\big( \mbf r, \mbf v \big)$ represents the dynamics, assumed as polynomial, or as the polynomial approximation of any selected application. 

The evaluation of the optimal control, at this step, follows classic derivations \cite{bryson2018applied}. The solution method consists into defining an auxiliary function, the control Hamiltonian:
\begin{align}
    \mathcal H(\mbf r,\mbf v, \tau, \hat{\boldsymbol \alpha} ) &= \mathcal P + \boldsymbol \lambda^T  \mbf p \big( \mbf x, \tau, \hat{\boldsymbol \alpha} \big) \\  
    &= \dfrac{1}{2}\mbf \tau^2 +  \boldsymbol \lambda_{\mbf r}^T \mbf v + \boldsymbol \lambda_{\mbf v}^T \left( \mbf f\big( \mbf r, \mbf v \big) + \tau \hat{\boldsymbol \alpha}  \right)
\end{align}
where vector $\boldsymbol \lambda$ represents the costates, with dimension equal to the number of states, divided into costates relative to position, $\boldsymbol \lambda_{\mbf r}$, and costates relative to the velocity, $\boldsymbol \lambda_{\mbf v}$. The equations of motion are obtained by applying the Hamilton's equations to the Hamiltonian, where optimal control is found through minimization.
\begin{equation}
    \begin{cases}
      \dot{\mbf r} = \dfrac{\partial \mathcal H}{\partial \boldsymbol \lambda_{\mbf r}} = \mbf v \\ 
      \dot{\mbf v}  = \dfrac{\partial \mathcal H}{\partial \boldsymbol \lambda_{\mbf v}} = \mbf f\big( \mbf r, \mbf v \big) + \tau \hat{\boldsymbol \alpha}  \\
      \dot{\boldsymbol \lambda_{\mbf r}}  = -\dfrac{\partial \mathcal H}{\partial \mbf r} = - \left[ \dfrac{\partial \mbf f}{\partial \mbf r} \right]^T \boldsymbol \lambda_{\mbf v}\\
      \dot{\boldsymbol \lambda_{\mbf v}}  = -\dfrac{\partial \mathcal H}{\partial \mbf v} = - \boldsymbol \lambda_{\mbf r} - \left[ \dfrac{\partial \mbf f}{\partial \mbf v} \right]^T \boldsymbol \lambda_{\mbf v}
    \end{cases}
\end{equation}
This system represents the augmented dynamics, where the time behaviour of the costates has been stacked onto the original set of ODEs. In order to find the optimal control, the terms in the Hamiltonian that depend only on $\tau$ and $\hat{ \boldsymbol \alpha}$ are isolated:
\begin{align}
    \mathcal H (\tau, \hat{\boldsymbol \alpha} ) = \dfrac{1}{2}\mbf \tau^2 +\boldsymbol \lambda_{\mbf v}^T  \tau \hat{\boldsymbol \alpha}  
\end{align}
By inspection, it is clear that the direction of trust that minimizes $\mathcal H$ is
\begin{equation}
    \hat{\boldsymbol \alpha}^* = - \dfrac{ \lambda_{\mbf v}}{ ||\lambda_{\mbf v}||}
\end{equation}
where starred variables indicate optimal values. By substituting the optimal direction back into the Hamiltonian 
\begin{equation}
    \mathcal H  = \dfrac{1}{2}\mbf \tau^2  - \tau ||\boldsymbol  \lambda_{\mbf v}||
\end{equation}
the necessary condition for the evaluation of $\tau$ is
\begin{equation}
    \dfrac{\partial \mathcal H}{\partial \tau} = 0 \quad \Rightarrow \quad \tau - ||\boldsymbol \lambda_{\mbf v}|| = 0 \quad \Rightarrow \quad \tau^* = ||\boldsymbol \lambda_{\mbf v}||
\end{equation}
Therefore, the optimal control law is explicit in $\mbf x$ and $\boldsymbol \lambda$
\begin{equation}
    \tau^* \hat{\boldsymbol \alpha}^* = ||\boldsymbol \lambda_{\mbf v}|| \left( - \dfrac{ \lambda_{\mbf v}}{ ||\lambda_{\mbf v}||} \right) = -\boldsymbol \lambda_{\mbf v} 
\end{equation}
and it can be substituted into the equations of motion:
\begin{equation}\label{sys}
    \begin{cases}
      \dot{\mbf r} = \dfrac{\partial \mathcal H}{\partial \boldsymbol \lambda_{\mbf r}} = \mbf v \\ 
      \dot{\mbf v}  = \dfrac{\partial \mathcal H}{\partial \boldsymbol \lambda_{\mbf v}} = \mbf f\big( \mbf r, \mbf v \big) - \boldsymbol \lambda_{\mbf v}   \\
      \dot{\boldsymbol \lambda_{\mbf r}}  = -\dfrac{\partial \mathcal H}{\partial \mbf r} = - \left[ \dfrac{\partial \mbf f}{\partial \mbf r} \right]^T \boldsymbol \lambda_{\mbf v}\\
      \dot{\boldsymbol \lambda_{\mbf v}}  = -\dfrac{\partial \mathcal H}{\partial \mbf v} = - \boldsymbol \lambda_{\mbf r} - \left[ \dfrac{\partial \mbf f}{\partial \mbf v} \right]^T \boldsymbol \lambda_{\mbf v}
    \end{cases}
\end{equation}
The optimal control has been resolved and the problem is now a TPBVP where the two initial states are known, $\mbf x(t_0) = [\mbf r(t_0) \quad \mbf v(t_0)]^T = \mbf x_0$ and $\mbf x(t_f) = [\mbf r(t_f) \quad \mbf v(t_f)]^T = \mbf x_f$, with the unknown of the initial costates $\boldsymbol \lambda(t_0) = [\boldsymbol \lambda_{\mbf r}(t_0) \quad \boldsymbol \lambda_{\mbf v}(t_0)]^T $, that propagated forward in time give exactly $\mbf x_f$. 

%%%%%%%%%%%%%%%%%%%%%%%%%%%%%%%%%%%%

\subsection{The Koopman Solution}
The equations of motion expressed in system \eqref{sys} can be solved using the KO approach. In the KO methodology, the solution of the set of ODEs is expressed as a linear combination of well-defined basis functions, where the evaluation of the Koopman matrix describes their time evolution and rate of change. Therefore, considering the augmented state composed of the original state variables with the addition of the costates, the solution of the state can be obtained as described by Eq.~\eqref{final}. Let us highlight the dependencies on time, state, and costate variables;
\begin{equation}\label{solution}
    \mbf x(\mbf x(t_0), \boldsymbol \lambda(t_0), t)
    = \mbf H \mbf V^{-1} \exp (\boldsymbol \Lambda t) \mbf V \mathcal L (\mbf x(t_0), \boldsymbol \lambda (t_0))
\end{equation}
where the matrix of observables $\mbf H$ has been calculated using the Galerkin method, after selecting the identity observables for the states. Equation \eqref{solution} is a polynomial map, function of $\mbf x_0$, $\boldsymbol \lambda_0$, and $t$, created by the basis functions, for each component of the state. In the TPBVP, we are interested in the values of the state at the given final time $t_f$, thus, the time dependency on the polynomials can be removed:
\begin{equation}
    \mbf x_f(\mbf x(t_0), \boldsymbol \lambda(t_0))
    = \mbf H \mbf V^{-1} \exp (\boldsymbol \Lambda t_f) \mbf V \mathcal L (\mbf x(t_0), \boldsymbol \lambda(t_0)) = \mathcal M_{t_0 \rightarrow t_f} (\mbf x(t_0), \boldsymbol \lambda(t_0)) \label{map}
\end{equation}
where $\mathcal M_{t_0 \rightarrow t_f}$ represents the polynomial map that connects any initial condition, $[\mbf x(t_0) \quad \boldsymbol \lambda(t_0) ]^T$, to the final state at the given time, $[\mbf x(t_f) \quad \boldsymbol \lambda(t_f) ]^T$. Therefore, the KO solution can be interpreted as a polynomial transition map up to the selected final time $t_f$. 

%%%%%%%%%%%%%%%%%%%%%%%%%%%%%%%%%%%%

\subsection{Control Solution Using Map Inversion}
The aim of the procedure is to find the value of costate at time $t_0$ that drives the state form $\mbf x_0$ to the desired final values $\mbf x_f$ at time $t_f$. However, the current KO solution, Eq. \eqref{map}, has $\boldsymbol \lambda (t_0)$ as one of the variables of the polynomials. Therefore, this map needs to be inverted, such that the initial conditions become a function of the final variables, and not vice-versa. In order for a map to be invertible, similarly to matrices, it needs to be square and full rank. The polynomials in $\mathcal M$ map $\mathbb{R}^{m} \rightarrow \mathbb{R}^{2m}$, where $m$ is the number of states. This issue is overcome by augmenting the solution and introducing the map corresponding to the identity function $\mbf x(t_0) = \mathcal I (\mbf x(t_0))$ \cite{valli2013nonlinear,di2014high}. This identity map of the state can be written highlighting the dependence on the costate variables, where all monomials connected to the costates are zero; $\mbf x(t_0) = \mathcal I_{\mbf x} (\mbf x(t_0), \boldsymbol \lambda (t_0))$. Thus, the augmented equation can be written as
\begin{equation}
    \begin{bmatrix}
    \mbf x(t_f) \\ \mbf x(t_0)
    \end{bmatrix}
    (\mbf x(t_0), \boldsymbol \lambda(t_0)) = 
    \begin{bmatrix}
    \mathcal M  \\ \mathcal I_{\mbf x}
    \end{bmatrix}
    (\mbf x(t_0), \boldsymbol \lambda(t_0))
\end{equation}
such that the augmented polynomial map meets the requirements to be inverted. Using inversion techniques for high order polynomials \cite{Berz1999ModernMM}, the map is inverted to obtain
\begin{equation}
    \begin{bmatrix}
    \boldsymbol \lambda(t_0) \\ \mbf x(t_0)   
    \end{bmatrix}
    (\mbf x(t_f), \mbf x(t_0)) = 
    \begin{bmatrix}
    \mathcal M  \\ \mathcal I_{\mbf x}
    \end{bmatrix}^{-1}
    (\mbf x(t_f), \mbf x(t_0))
    = \mathcal W (\mbf x(t_f), \mbf x(t_0))
\end{equation}
The inverted map $\mathcal W$ is a function of two variables: the state at the initial time $\mbf x(t_0)$, and the values of the states at final time $\mbf x_f(t_f)$. The initial costates are found by considering the first half of the inverted $\mathcal W$ map, since the first $m$ polynomials are related to the control:
\begin{equation}
    \boldsymbol \lambda(t_0) 
    = \mathcal W_{\boldsymbol \lambda} (\mbf x(t_f), \mbf x(t_0)) \label{eq:findcostate}
\end{equation}
where $\mathcal W_{\boldsymbol \lambda}$ indicates the map with only the costate polynomials selected. The optimal numerical values are obtained through polynomial evaluation
\begin{equation}
    \boldsymbol \lambda_0
    = \mathcal W_{\boldsymbol \lambda} (\mbf x_f, \mbf x_0) \label{eva}
\end{equation}
where the substitutions $\mbf x(t_f) = \mbf x_f$ and $\mbf x(t_0) = \mbf x_0$ guarantee to satisfy the boundary conditions of the TPBVP. The control optimization is concluded and the time evolution evolution of both the states, and the costates, can be found by solving system \eqref{sys} with initial conditions $\mbf x(t_0) = \mbf x_0$ and $\boldsymbol \lambda(t_0) = \boldsymbol \lambda_0$.

It is important to underline that contrary to other techniques based on map inversion \cite{valli2013nonlinear,di2014high}, the KO solution does not require the evaluation of any nominal control to perform accurately. This aspect is due to the global nature of the KO approximation, and to the fact that the KO map of the solution is able to work directly with state variables over deviations. Looking at Eq.~\eqref{eq:findcostate}, the inverted map creates a function of the optimal costates on the boundary conditions. Therefore, if we seek a different solution where either the initial condition and/or the desired final state are changed, the optimal control can be rapidly obtained by evaluating the map as in Eq.~\eqref{eva}. The Koopman approach has completely described the influence of the boundary conditions on the optimal initial costate, such that once the Koopman solution is obtained, any change in the application can be resolved by looking directly at the state transition polynomial maps.

%%%%%%%%%%%%%%%%%%%%%%%%%%%%%%%%%%%%
%%%%%%%%%%%%%%%%%%%%%%%%%%%%%%%%%%%%
%%%%%%%%%%%%%%%%%%%%%%%%%%%%%%%%%%%%

\section{Numerical Applications}
The proposed control techniques is applied to three different numerical applications. First, the state of the Duffing oscillator is driven to rest from multiple configurations, displaying the control effort and the relationship with the maneuver time. Afterwards, the Clohessy-Wiltshire (CW) problem is presented and the equations of motion are derived. This application is used to prove the efficiency of the new technique for rendezvous problems. Lastly, the assumptions that lead to the classical linear CW equations are relaxed and a new, high-order, system that better describes the relative motion between two satellites is presented, where nonlinearites play an important role in increasing the accuracy of the relative motion. 

%%%%%%%%%%%%%%%%%%%%%%%%%%%

\subsection{The Duffing Oscillator}
The application to the Duffing Oscillator system is now proposed, where the KO inverse control technique is applied to control the state to a set of desired values for different time requirements. The Duffing oscillator describes the oscillations of a mass $M$ attached to a nonlinear spring,  with stiffness constant $k$, and a damper. The scalar system can therefore be described by the following equations of motion, where $q$ represents the position, while $p$ is the velocity;
\begin{subequations} 
\begin{align} \
    \dot q &= \dfrac{p}{M}\\
    \dot p &= -k q - k a^2\epsilon q^3
\end{align}
\end{subequations}
where $a$ is a unit transformation constant and $\epsilon$ is a small parameter. In this scenario, the ODE parameters are defined as $a = 1$, $M = 1$, $k = 1$, and $\epsilon = 0.001$. The system is at equilibrium at rest, for null values of position and velocity, $(0,0)$. When perturbed, the mass starts to oscillate in a periodic motion. After adding control to the oscillator, the controlled system, in the form of Eqs. \eqref{sys}, has the following form:
\begin{subequations} 
\begin{align} \
    \dot q &= p\\
    \dot p &= - q - \epsilon q^3 - \lambda_p\\
    \dot \lambda_q &= 3\epsilon q^2 \lambda_p + \lambda_p\\
    \dot \lambda_p &= - \lambda_q
\end{align}
\end{subequations}

\begin{figure}[h]
\centering
\includegraphics[width=\textwidth]{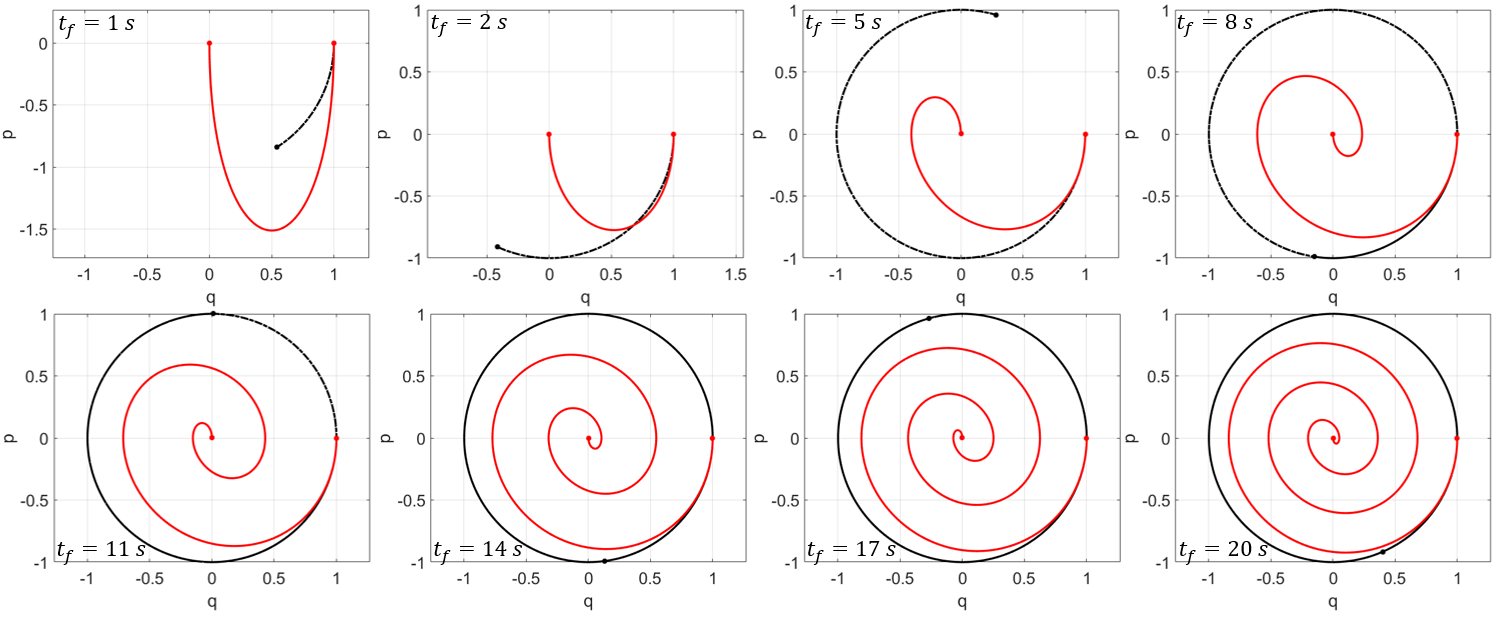}
\caption{Bringing the Duffing oscillator to its equilibrium point for different times.}
\label{fig:times}
\end{figure}
At first, let us assume that the mass is moved from its equilibrium position to $q_0 = 1$, and given no initial velocity, $p_0 = 0$. The goal of this application is to bring the mass back at rest in its equilibrium point given a fixed time length. Figure \ref{fig:times} shows the system brought to the origin for different amounts of times. The uncontrolled systems has the mass oscillating between -1 and 1 indefinitely, displayed in black. The energy optimal approach has the state of the system rotating in a spiral around the equilibrium point, at a rate that gets slower the larger the final time $t_f$. 
\begin{figure}[h]
\centering
\includegraphics[width=.7\textwidth]{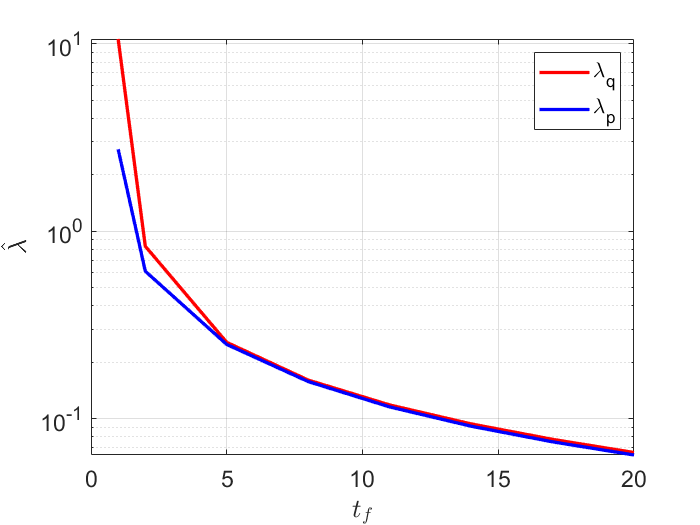}
\caption{Control magnitude requested for each time length.}
\label{fig:mag}
\end{figure}
Thus, the required magnitude of control effort decreases as the system is allowed more time to reach stasis. Figure \ref{fig:mag} reports the amount of control applied to the system for different times, both for position and velocity, evaluated as
\begin{equation}
    \hat{\lambda_i} = \dfrac{1}{t_f}\sqrt{\int_0^{t_f} \lambda_i^2 dt}
\end{equation}
with $i=q,p$. The control effort monotonically decreases as time increases, with the costate for position and velocity overlapping for large amount of times. 

\begin{figure}[h]
\centering
\includegraphics[width=.75\textwidth]{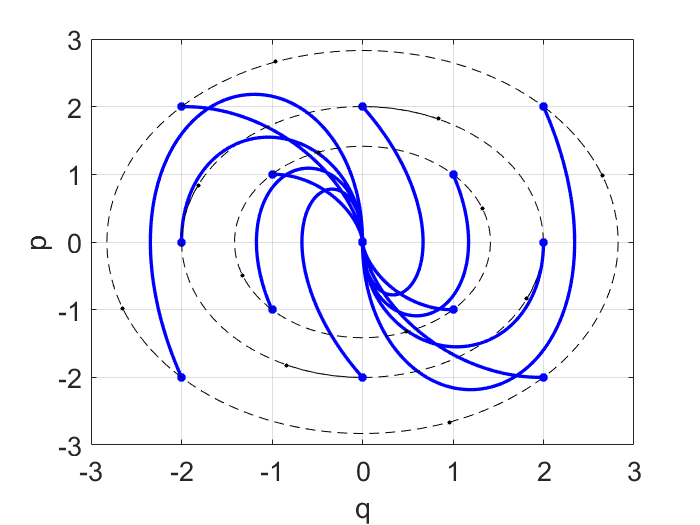}
\caption{Controlled pathways to equilibrium for $t_f = 2$ seconds.}
\label{fig:grid}
\end{figure}
A second analysis is performed with a fixed required time $t_f = 2$ seconds, where the initial condition of the state is picked following a grid, and the system is lead to stasis. Figure \ref{fig:grid} shows the path to equilibrium from the 12 different initial conditions, and the relative uncontrolled behaviour in the background. Therefore, control has been achieved given different states and time constraints. 

%%%%%%%%%%%%%%%%%%%%%%%%%%%

\subsection{The Clohessy-Wiltshire problem}
The Clohessy-Wiltshire (CW) equations represent a linear approximation of the motion between two, relatively close, satellites in orbit \cite{curtis2013orbital}. The problem can be formulated using different set of coordinates \cite{kasdin2005canonical}: the approach of deriving the Lagrangian for the relative motion in Cartesian coordinates is here chosen. Then, the Hamiltonian of the relative system can be evaluated using a Legendre transformation, such that it is possible to divide the linear part from the nonlinear contribution. 

\begin{figure}[h]
\centering
\includegraphics[width=.65\textwidth]{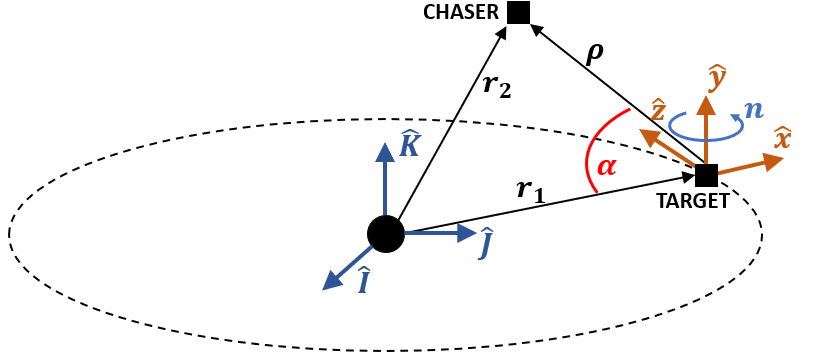}
\caption{Relative motion rotating frame.}
\label{fig:frame}
\end{figure}

The Cartesian coordinates are the most convenient when allowing control and simulation techniques. Let us consider the rotating Cartesian Euler-Hill system, where the origin of the frame is set on a circular orbit of radius $a$ about the primary celestial body. Therefore, this frame is rotating with mean motion $n = \sqrt{\dfrac{\mu}{a^3}}$, where $\mu$ is the gravitational parameter. The classic Clohessy-Wiltshire equations are obtained with the assumption of relative motion of two spacecrafts with respect to a circular orbit. Let us define with $\mbf r_1$ the position vector of the target spacecraft, and with $\mbf r_2$ the position vector of the chaser, as shown in Figure \ref{fig:frame}, with respect to the $\mathcal{I} = \{\hat I,\hat J, \hat K\}$ reference frame. 
Thus, the relative position vector in the rotating frame $\mathcal{R}=\{\hat x,\hat y, \hat z\}$ is defined as $\mbf \rho = [x,y,z]^T$, and the angular velocity between the two frames $\mathcal{I}$ and $\mathcal{R}$ is $\boldsymbol \omega = [0,0,n]^T$. Under the circular orbit assumption, $||\mbf r_1|| = a$, and the component-wise relative velocity in the rotating frame can be calculated by
\begin{equation}
    \mbf v = \boldsymbol \omega \times \mbf r_1 + \boldsymbol \omega \times \boldsymbol \rho + \dfrac{d^{\mathcal{R}}}{dt} \boldsymbol \rho = 
    \begin{bmatrix}
    \dot x - n y \\ \dot y + n x + n a \\ \dot z
    \end{bmatrix}
\end{equation}
The equations of motion are obtained by applying the Euler-Lagrange equation to the Lagrangian of the system $\mbf L$. The Lagrangian is evaluated by subtracting the potential energy $\mbf U$ to the kinetic energy $\mbf T$:
\begin{equation}
    \mbf L = \mbf T - \mbf U
\end{equation}
Working with quantities per unit of mass, the kinetic energy is:
\begin{equation}
    \mbf T = \dfrac{1}{2}||\mbf v|| = \dfrac{1}{2}(\dot x +\dot y + \dot z)^2 + n(x\dot y - y \dot x + a\dot y ) +\dfrac{1}{2}n^2y^2 +\dfrac{1}{2}n^2x^2 +n^2 a x +\dfrac{1}{2}n^2 a^2 
\end{equation}
The potential energy of the chaser is the one of a spacecraft attracted by a spherical body. This gravitational potential is expanded using Legendre polynomials following the developments from \cite{kasdin2005canonical}:
\begin{equation}
    \mbf U = -\dfrac{\mu}{||\mbf r_2||} = -\dfrac{\mu}{||\mbf r_1 + \boldsymbol \rho||} = 
    -\dfrac{\mu}{a\sqrt{1+2 \dfrac{\mbf r_1 \cdot \boldsymbol \rho}{a^2} + \dfrac{\rho^2}{a^2} }} =
     -\dfrac{\mu}{a}\sum_{i=0}^{\infty} \mathcal P_k(\cos \alpha) \left(\dfrac{\rho}{a} \right)^k
\end{equation}
where $\rho = ||\boldsymbol \rho||$, and $\mathcal P_k(\cos \alpha)$ are the Legendre polynomials. Looking back at Figure \ref{fig:frame}, angle $\alpha$ is the angle between the position of the target and the relative position of the chaser. Therefore, the cosine can be rewritten as 
\begin{equation}
    \cos \alpha = -\dfrac{\mbf r_1 \cdot \boldsymbol \rho}{a\rho} = -\dfrac{x}{\rho}
\end{equation}
Substituting in the Lagrangian, the formulation becomes
\begin{equation}\label{eq:lag}
    \mbf L = \dfrac{1}{2}(\dot x +\dot y + \dot z)^2 + n(x\dot y - y \dot x + a\dot y ) +\dfrac{1}{2}n^2y^2 +\dfrac{1}{2}n^2x^2 +n^2 a x +\dfrac{1}{2}n^2 a^2 + n^2a^2 \sum_{k=0}^{\infty} \mathcal P_k \left(-\dfrac{x}{\rho}\right) \left(\dfrac{\rho}{a} \right)^k
\end{equation}
which, in practice, needs to be approximated by selecting the maximum order of the terms in the summation. Therefore, it is convenient to define the set of polynomials $\mathcal Q_k$ as
\begin{equation}
    \mathcal Q_k = n^2a^2 \mathcal P_k\left(-\dfrac{x}{\rho}\right) \left(\dfrac{\rho}{a} \right)^k  \quad \forall \quad k=0,\dots,\infty
\end{equation}
such that it eases the contribution of each singular term added to the Lagrangian by the potential for each order. Table \ref{tab:table1} shows the contribution of the first 6 polynomials.
\begin{table}
\caption{\label{tab:table1} Potential polynomials for different orders.}
\centering
\begin{tabular}{c|c}
\hline
$k$ & $\mathcal Q_k$ \\
\hline
0 & $n^2a^2$ \\
1 & $-n^2 a x$ \\
2 &  $\dfrac{3}{2} n^2 x^2 -  \dfrac{1}{2}\rho^2 a^2$\\
3 &  $- \dfrac{5}{2} \dfrac{n^2}{a}x^3 + \dfrac{3}{2} \dfrac{\rho^2n^2}{a}x $\\
4 &  $  \dfrac{35}{8} \dfrac{n^2}{a^2}x^4 - \dfrac{30}{8} \dfrac{\rho^2 n^2}{a^2}x^2 + 3 \dfrac{\rho^4n^2}{a^2}$\\
5 &  $ -\dfrac{63}{8} \dfrac{n^2}{a^3}x^5 + \dfrac{70}{8} \dfrac{\rho^2n^2}{a^3}x^3 - \dfrac{15}{8} \dfrac{\rho^4n^2}{a^3}x   $\\
$\dots$ & $\dots$\\
\hline
\end{tabular}
\end{table}
In the Clohessy-Wiltshire equations for relative motions, only small deviations from the reference orbit are considered, where only the first three polynomials of the potential expansion are considered. Thus, substituting the polynomials with $\rho = \sqrt{x^2 + y^2 + z^2}$, we obtain the low order Lagrangian
\begin{equation}
    \mbf L_{LIN} = \dfrac{1}{2}(\dot x +\dot y + \dot z)^2 + n(x\dot y - y \dot x + a\dot y )
    +\dfrac{3}{2}n^2 a^2 +\dfrac{3}{2}n^2x^2 - \dfrac{1}{2}n^2z^2
\end{equation}
where the subscript ``$LIN$" indicates that this Lagrangian provides the linear equations of motion. In fact, the classic linear CW system of ODE is obtained by applying the Euler-Lagrange equation
\begin{equation}
    \dfrac{\partial \mbf L}{\partial x} - \dfrac{d}{d t} \left(\dfrac{\partial \mbf L}{\partial \dot x}\right) = 0
\end{equation}
repeated for variables $y$ and $z$. The usual equations of motion are then obtained:
\begin{subequations} \label{linsy}
\begin{align}
    \ddot x &= 2n\dot y + 3n^2 x \\
    \ddot y &= -2n\dot x \\
    \ddot z &= n^2 z
\end{align}
\end{subequations}
However, it is common procedure to normalize these equations and work with a simplified notation. Therefore, all distances get normalized by the semi-major axis $a$, which corresponds to a normalization of the rates of a factor $n$. Under this transformation, the new variables are dimensionless and the Lagrangian can be rewritten as 
\begin{equation}
    \tilde{\mbf L}_{LIN} = \dfrac{1}{2}(\dot x +\dot y + \dot z)^2 + x\dot y - y \dot x + \dot y
    +\dfrac{3}{2} +\dfrac{3}{2}x^2 - \dfrac{1}{2}z^2
\end{equation}
which gives the normalized system of ODEs:
\begin{subequations} 
\begin{align}
    \ddot x &= 2\dot y + 3 x \\
    \ddot y &= -2\dot x \\
    \ddot z &= z
\end{align}
\end{subequations}
It can be noted that the evolution of the system in the $z$ direction is completely decoupled from the $(x,y)$ plane motion and it could be solved by considering the contribution separately. In the KO solution, the system gets further scaled by a factor of 0.05 to best fit the evaluation of each component of the Koopman matrix with the Galerkin method.   

\begin{figure}[h]
\centering
\includegraphics[width=.8\textwidth]{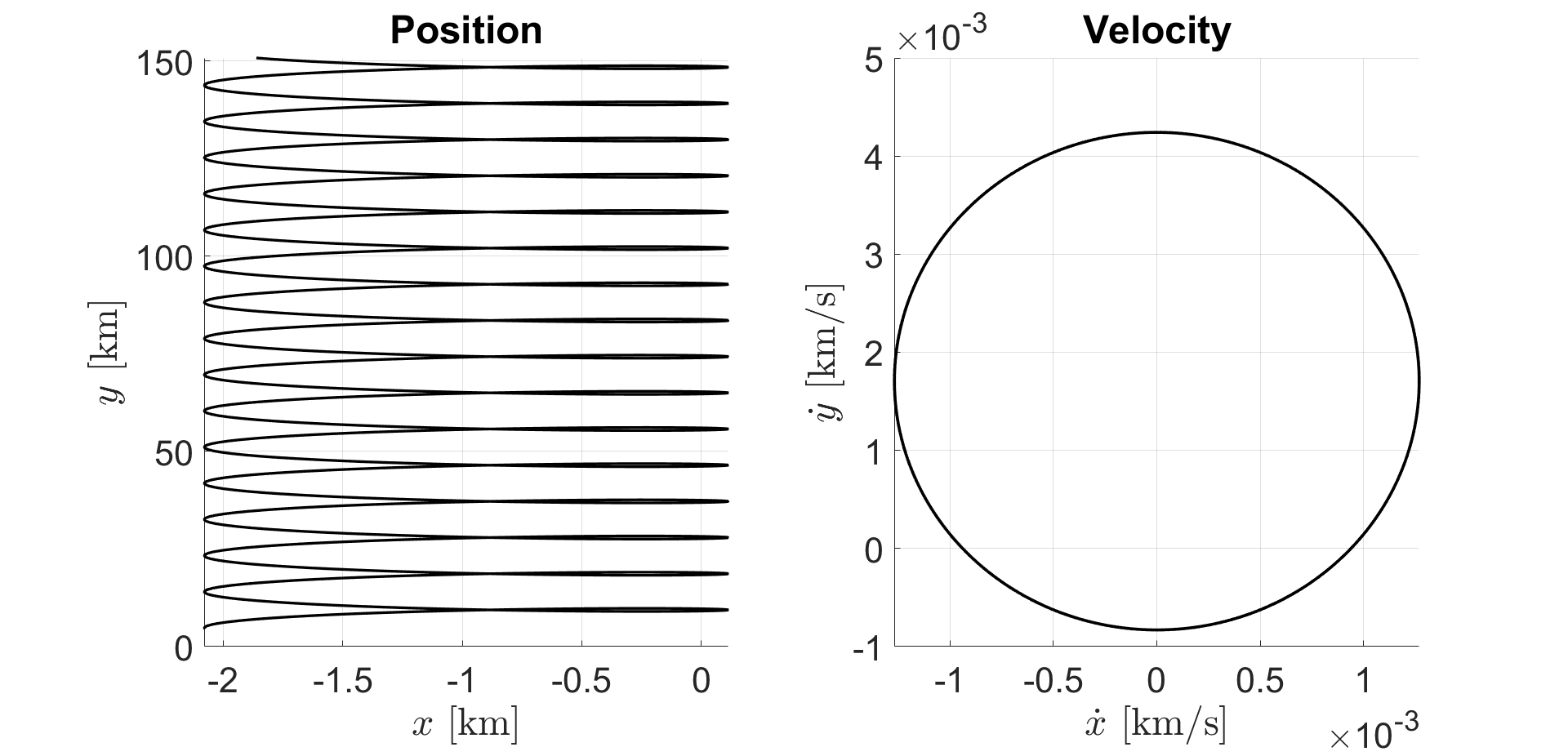}
\caption{Relative uncontrolled motion of the chaser with respect to the target for one day.}
\label{fig:uncon}
\end{figure}

\begin{figure}[h]
\centering
\includegraphics[width=.9\textwidth]{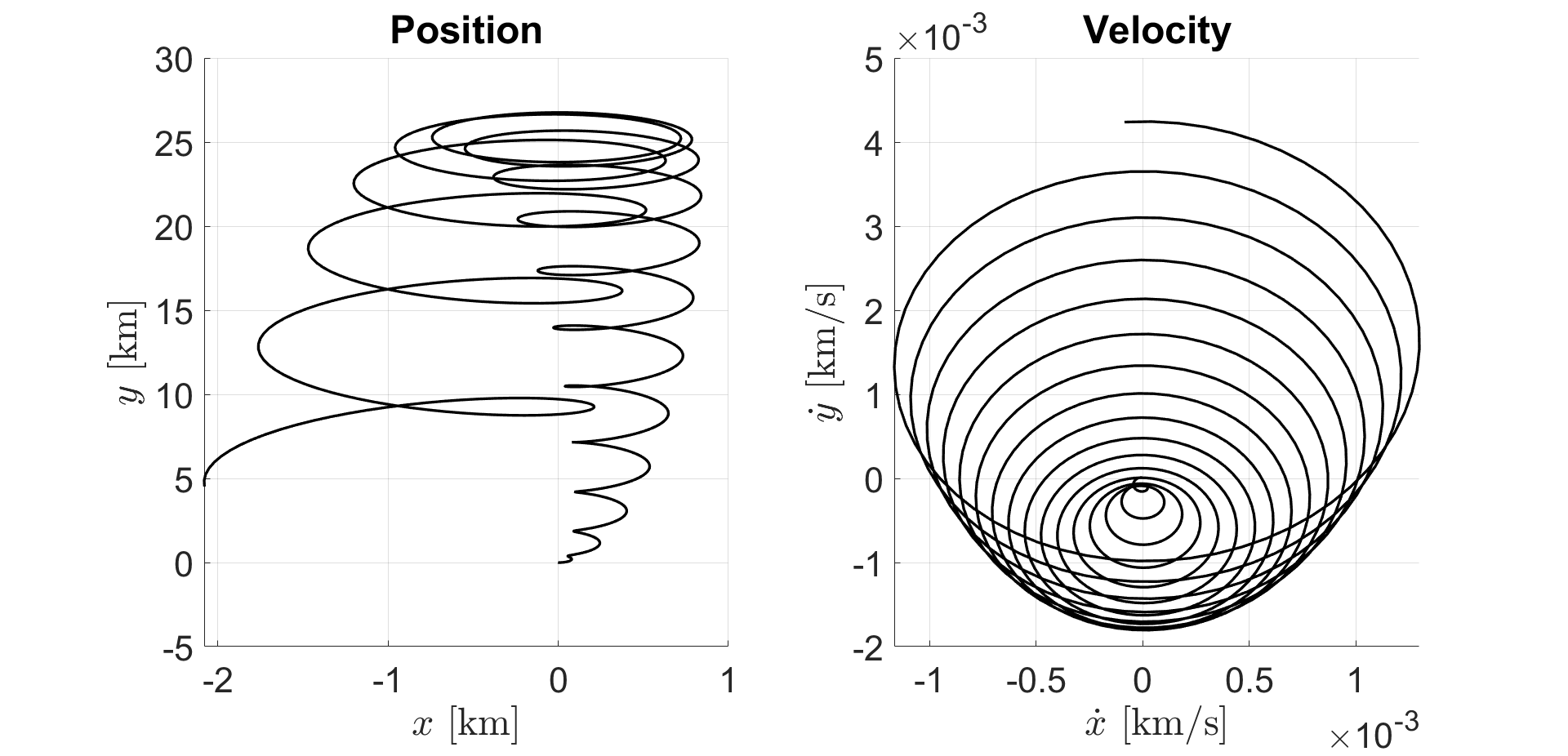}
\caption{Rendezvous of the chaser to the target in one day.}
\label{fig:con}
\end{figure}

Therefore, consider a target spacecraft orbiting around earth with a semi-major axis $a = 6678000$ meters. In this application, a chaser satellite aims to rendezvous with the target in one day. The relative (unscaled) initial state of the chaser with respect to the target is 
\begin{subequations} 
\begin{align}
    \mbf r_0 &=
    \begin{bmatrix}
    -2.0772 & 4.5157 & 0 
    \end{bmatrix} \text{km} \label{ic1}\\
    \mbf v_0 &=
    \begin{bmatrix}
    -8.6074\times 10^{-5} & 4.2376\times 10^{-3} & 0 
    \end{bmatrix} \text{km/s} \label{ic2}
\end{align}
\end{subequations}
which means that the chaser is getting farther away from the target, and that it needs to slow down and invert its relative direction of motion. The uncontrolled motion for one day is displayed in Figure \ref{fig:uncon}, where the left graph shows the relative position of the chaser with respect to the target in the $(x,y)$ plane, while the graph on the right portraits the relative velocity profile. 
The uncontrolled chaser keeps getting away from the target in a periodic motion. Thus, control is added to the system in an energy optimal manner, such that the controlled dynamics are written in the formulation provided by Equations \eqref{sys}: 
\begin{subequations} \label{aug}
\begin{align}
    \dot x &= v_x\\
    \dot y &= v_y \\
    \dot z &= v_z \\
    \dot v_x &= 3n^2 x + 2 n v_y - \lambda_{v_x}\\
    \dot v_y &= -2n v_x - \lambda_{v_y}\\
    \dot v_z &= -n^2 z - \lambda_{v_z}\\
    \dot \lambda_x &= -3n^2\lambda_{v_x}\\
    \dot \lambda_y &= 0 \\
    \dot \lambda_z &= n^2 \lambda_{v_z} \\
    \dot \lambda_{v_x} &= -\lambda_{x}+2n \lambda_{v_y}\\
    \dot \lambda_{v_y} &= -\lambda_{y}-2n \lambda_{v_x}\\
    \dot \lambda_{v_z} &= -\lambda_{z}
\end{align}
\end{subequations} 
After solving the system in the KO framework, and inverting the polynomial map of the augmented state (order three selected), the values of the initial costates that bring the system to rendezvous are
\begin{align}
    \boldsymbol \lambda_{\mbf r,0} &=
    \begin{bmatrix}
    -4.3659\times 10^{-11} & 1.6400\times 10^{-13} & 0 
    \end{bmatrix}^T \label{lam1}\\
    \boldsymbol \lambda_{\mbf v,0} &=
    \begin{bmatrix}
    -1.0017\times 10^{-9} & -1.5900\times 10^{-8} & 0 
    \end{bmatrix}^T \label{lam2}
\end{align}
The controlled pathway of the chaser can be appreciated in Figure \ref{fig:con}, where it can be noted that the spacecraft slows down and changes its direction to approach the origin with smaller oscillations each time. The velocity profile, the right plot of the figure, shows how the chaser decreases its velocity in a spiral, reaching stasis at the desired final time. 
The rendezvous is achieved and the control technique has been proven to work successfully. 

In this specific case that system \eqref{aug} is linear, it is possible to find the analytical solution of the inverse control problem. Indeed, let us consider the state transition matrix (STM) $\boldsymbol \Phi $ such that
\begin{equation}
    \begin{bmatrix}
    \mbf x_f \\ \boldsymbol \lambda_f
    \end{bmatrix}
    = \begin{bmatrix}
    \boldsymbol \Phi_{11} & \boldsymbol \Phi_{12}\\ \boldsymbol \Phi_{21} & \boldsymbol \Phi_{22}
    \end{bmatrix}
    \begin{bmatrix}
    \mbf x_0 \\ \boldsymbol \lambda_0
    \end{bmatrix}
\end{equation}
where the STM is calculated as $\boldsymbol \Phi = \exp{(\mbf A t_f)} $, with matrix $\mbf A$ obtained by form the set of Equations \eqref{aug}. Therefore, for the linear case, it is possible to solve for the analytical solution of the initial costate as
\begin{equation}
    \boldsymbol \lambda_0 = \boldsymbol \Phi_{12}^{-1}\left( \mbf x_f -  \boldsymbol \Phi_{11} \mbf x_0 \right) \label{anal}
\end{equation}
When comparing the KO solution, Equations \eqref{lam1} and \eqref{lam2}, with the analytical solution obtained through Equation \eqref{anal}, we obtain the following levels of relative errors, for the position $\epsilon_{REL}(\boldsymbol \lambda_{\mbf r,0})$, and velocity $\epsilon_{REL}(\boldsymbol \lambda_{\mbf v,0})$;
\begin{align}
    \epsilon_{REL}(\boldsymbol \lambda_{\mbf r,0}) &=
    \begin{bmatrix}
    0.021 & 0.001 & 0 
    \end{bmatrix}^T \%\\
    \epsilon_{REL}(\boldsymbol \lambda_{\mbf v,0}) &=
    \begin{bmatrix}
    0.840 & 0.033 & 0 
    \end{bmatrix}^T \%
\end{align}
proving a correct evaluation of the optimal costates. This comparison validates the reliability of the KO methodology, that can be now applied to the nonlinear CW dynamics. 

%%%%%%%%%%%%%%%%%%%%%%%%%%%

\subsection{The High-Order Clohessy-Wiltshire Problem}
The relative equations of motion studied so far consider only small deviations of the chaser from the reference orbit. However, when the magnitude of the relative position and velocity increases, the linear approximation of the Clohessy-Wiltshire equations is not sufficient to fully describe the system. Therefore, in this new application, the effects of the nonlinearities between chaser and target are taken into account. Thus, looking back at the definition of the Lagrangian from Equation \eqref{eq:lag}, the potential component can be truncated at orders higher than two (which gives the linear system). The resulting high-order Clohessy-Wiltshire system, obtained via the Euler-Lagrange equations, is fully coupled. 

Looking at Table \ref{tab:table1}, for any order $k$, the polynomials $\mathcal Q_k$ have the magnitude of the relative position vector, $\rho$, appearing only with even powers. Indeed, when looking at the $i$-th polynomial, the monomials of  $\mathcal Q_i$ are in the form
\begin{align}
  \mathcal Q_i &= q_1 \rho^i + q_2 \rho^{(i-2)}x^2 + q_3 \rho^{(i-4)}x^4 + \dots  + q_n \rho^0 x^i \quad \forall i=2,4,6,\dots \\
  \mathcal Q_i &= q_1 \rho^{(i-1)} + q_2 \rho^{(i-3)}x^3 + q_3 \rho^{(i-5)}x^5 + \dots  + q_n \rho^0 x^i \quad \forall i=3,5,7,\dots
\end{align}
which means that any $\mathcal Q_k$ is polynomial in the variables $x$, $y$, and $z$, since the square root of $\rho = \sqrt{x^2 + y^2 + z^2}$ is always cancelled by the even power of its exponent. Therefore, the resulting Lagrangian for any truncation order will always give polynomial equations of motion from the Euler-Lagrangian equations, which are suitable to the KO control technique. While the truncation order of the potential energy can be chosen arbitrary large, the effects of the nonlinear terms and the accuracy of the model reaches a plateau for large numbers. Indeed, looking again back at Table  \ref{tab:table1}, it is worth noticing the the magnitude of the coefficients of the monomials of $\mathcal Q_k$ decreases by a factor $a$ for each additional order. That is, for a large $k$, the contribution of the new polynomials become negligible when compared to the low order polynomials. The constant decrease of the magnitude of the coefficients is directly related by the size of the target's orbit. Since each high-order polynomial gets divided by the semi-major axis, the effect of nonlinearities is stronger for low altitude orbits. Conceptually, a stronger curvature of the orbit, closer to Earth, has stronger nonlinear effects when compared to far orbits, where the local linearization is more accurate. Therefore, the high-order Clohessy-Wiltshire equations present their advantages and improvements in the LEO region. 

\begin{figure}[h]
\centering
\includegraphics[width=.75\textwidth]{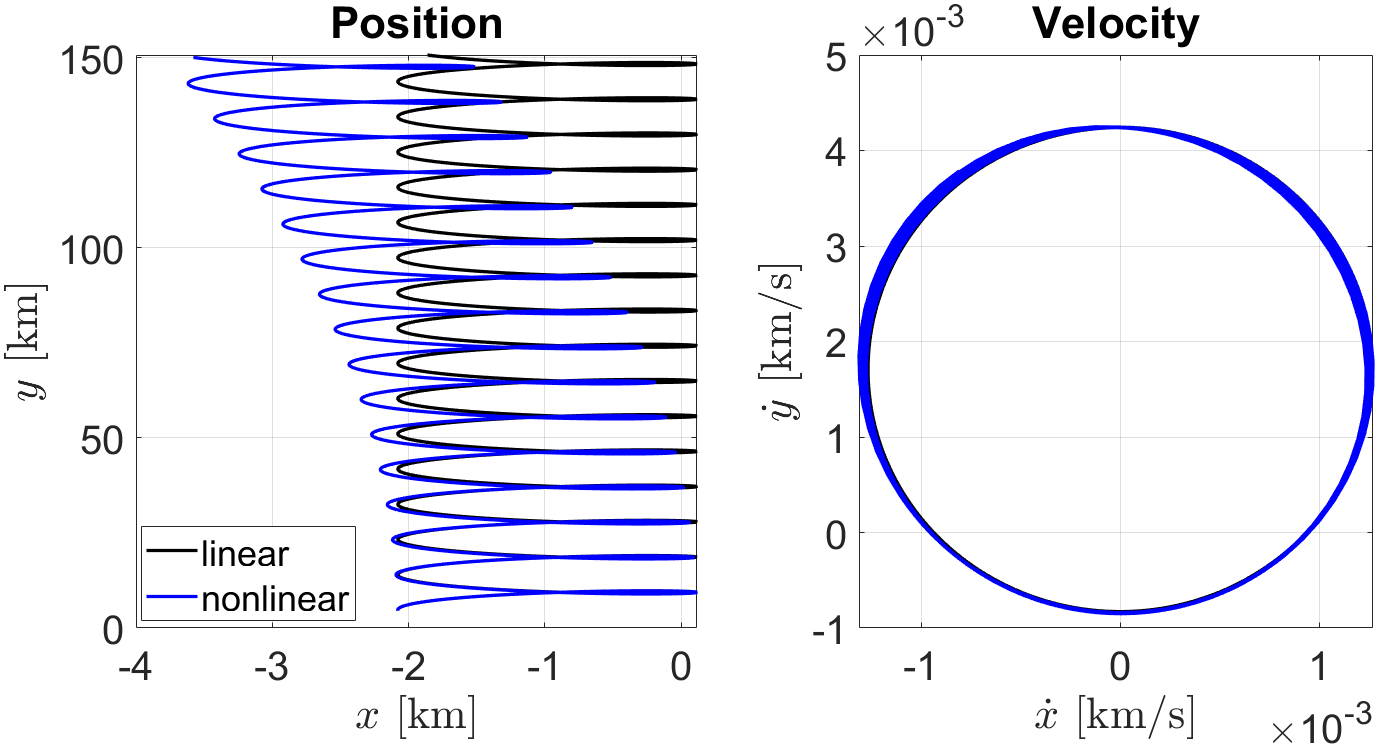}
\caption{Comparison between classic and high-order Clohessy-Wiltshire equations.}
\label{fig:vs}
\end{figure}
Let us consider the initial condition previously selected by Eq. \eqref{ic1} and Eq. \eqref{ic2}. For the simulations proposed in this paper, the target is orbiting in LEO at 307 km altitude. Figure \ref{fig:vs} shows the comparison between the classic (linear) Clohessy-Wiltshire propagation, in black lines, and the high-order Clohessy-Wiltshire propagation, in blue. The linear propagation is analogous at Figure \ref{fig:uncon}, which correspsonds to selecting truncation order two for the potential polynomials. The nonlinear (blue) lines are instead evaluated by doubling the truncation order to four. Looking at the position graph on the left, after an initial overlapping between the linear and nonlinear solution, the blue line starts to drift apart transitioning on the left, while the black line keeps a perfect periodic behaviour. Therefore, the high-order solution better approximates the actual relative position of the chaser, as shown in \cite{curtis2013orbital}. 

\begin{figure}[h]
\centering
\includegraphics[width=.8\textwidth]{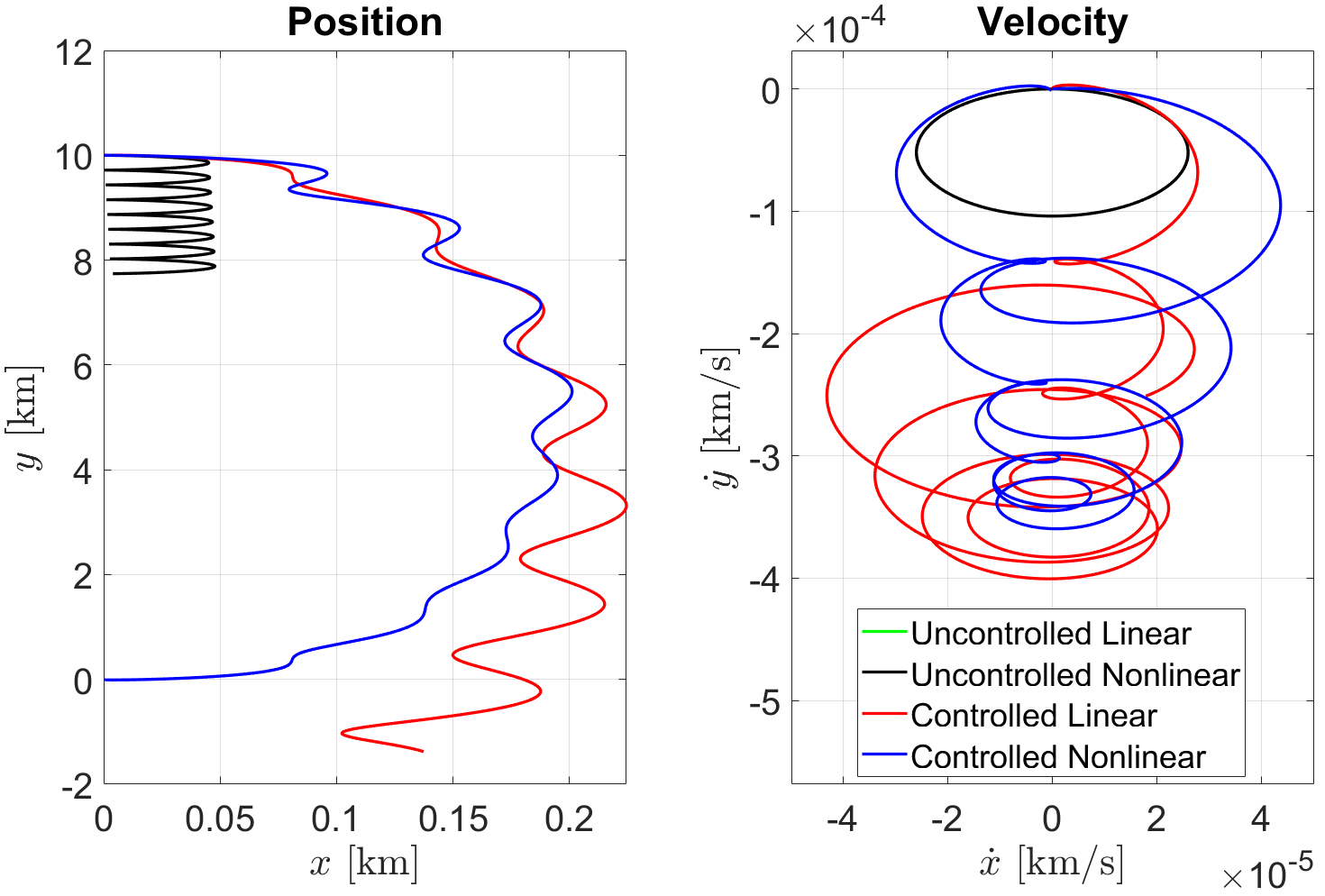}
\caption{Rendezvous of the chaser to the target in 12 hours under different models.}
\label{fig:noncon}
\end{figure}

The KO control technique is applied to reach rendezvous between chaser and target when the high-order dynamics is considered. In order to show the benefits of the new approximation, an initial condition of the type
\begin{subequations} 
\begin{align}
    \mbf r_0 &=
    \begin{bmatrix}
    0 & $y$ & 0 
    \end{bmatrix} \text{km}, \\
    \mbf v_0 &=
    \begin{bmatrix}
    0 & 0 & 0 
    \end{bmatrix} \text{km/s}
\end{align}
\end{subequations}
has been selected. This initial condition underlines the benefits of the nonlinear systems since, in the linear dynamics, the $y$ axis is a set of equilibrium points. Indeed, looking at system \eqref{linsy}, any point in the $y$ axis is stationary and the forward propagation of such initial condition leads to a constant vector. Physically, this means that if the chaser is exactly behind the target with null relative velocity, its relative position is fixed in time. Such approximation holds only when the curvature of the orbit is neglected (full linearization). Therefore, the aim of the presented application is to achieve rendezvous when the chaser is at $y=10$ km ahead of the target, given the constraint of performing the maneuver in 12 hours. Figure \ref{fig:noncon} shows the resulting pathways of the chaser subject to different dynamics. As stated above, the uncontrolled linear propagation results to a constant integration, that is a single point in the figure. Therefore, the uncontrolled linear solution, in green, is not visible since it is static. On the contrary, the uncontrolled nonlinear dynamics, in black, has been evaluated considering the first six potential polynomials. This line shows the actual behaviour related to the uncontrolled relative motion of the spacecraft, with the presence of the translational drift portrayed in Figure \ref{fig:vs}. The controlled pathways are reported in the figure with red and blue lines, for the linear and the high-order dynamics, respectively. Thus, the correct solution of the controlled chaser can be appreciated by analyzing the blue curves. As the chaser gains velocity by increasing its $x$ components, it is able to rapidly approach the target in the $y$ direction. After passing the half-point, the control starts to slow the chaser down, such that it reaches rendezvous with null velocity. The controlled linear pathway, in red, is reported for comparison purposes. The red curve shows the behaviour of the chaser that is propagated under linear dynamics, which fails to achieve rendezvous. Thus, the consideration of nonlinear terms is crucial for the high accuracy and precision evaluation of the costates of the system and the newly developed energy optimal control technique is able to evaluate such control. 

\begin{figure}[h]
\centering
\includegraphics[width=\textwidth]{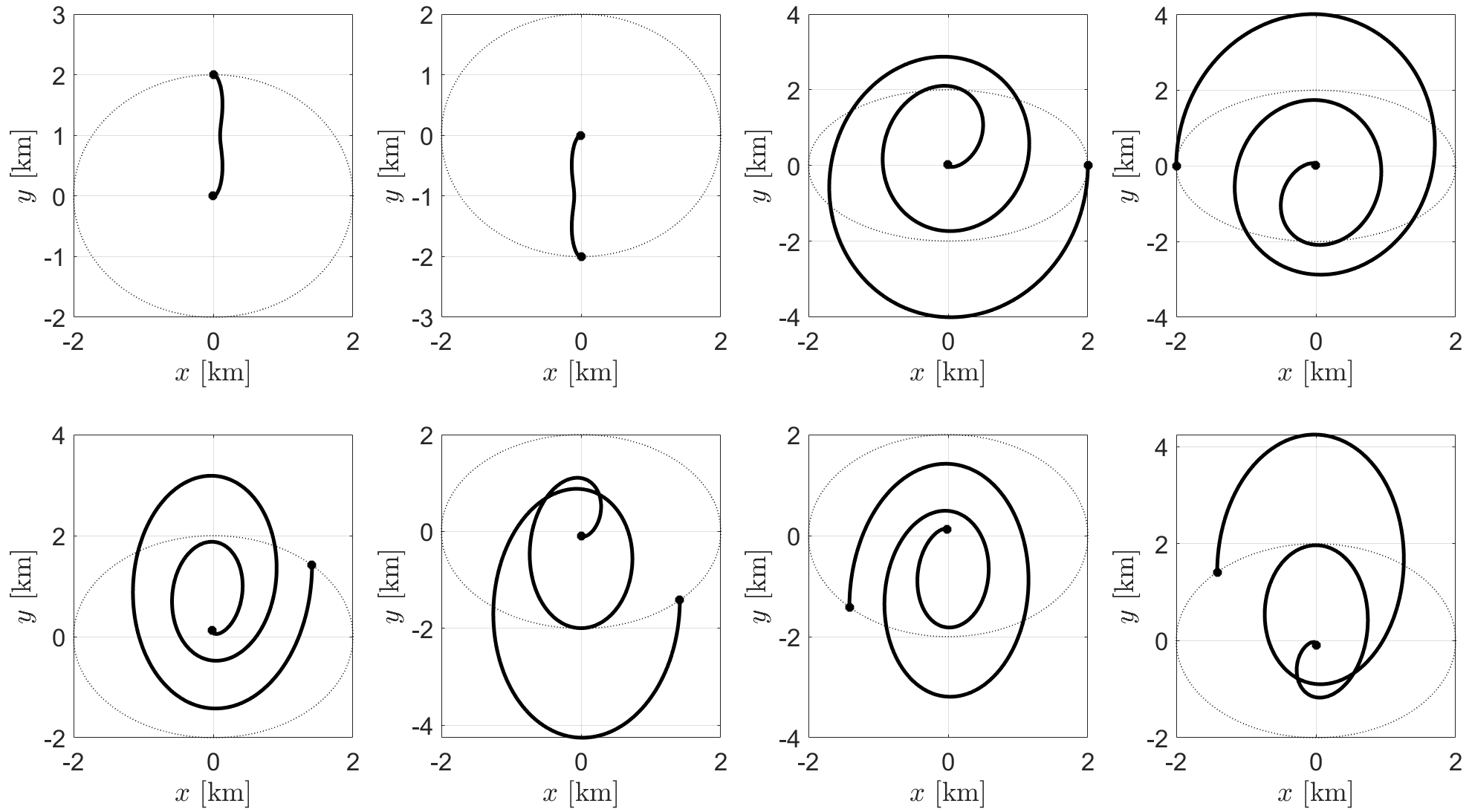}
\caption{Rendezvous of the chaser to the target in 4 hours given different relative positions in orbit.}
\label{fig:gridNON}
\end{figure}

The rendezvous application has been repeated with the chaser initial condition at different relative positions, in order to asses the capability of the KO technique to evaluate the optimal control to approach the target from every direction. As such, Figure \ref{fig:gridNON} reports the optimal rendezvous pathways of the chaser starting from a distance of two kilometers from the target (dotted circle), with initial condition separated every $\pi/4$. The figure reports solutions for rendezvous in exactly four hours. The values of the initial relative velocity $\dot y$ is evaluated considering the difference in the orbit's tangential velocity, where the reference orbit radius $a$ is increased by a factor $x$. This aspect best describes the relative initial condition between chaser and target orbiting in different orbits. Therefore, Figure \ref{fig:gridNON} shows that the KO energy optimal inverse control technique is able to calculate the exact costates for the chaser approaching from every direction. These solutions have been calculated using a high-order Lagrangian that considers the first six polynomials in the expansion of the potential energy. 

The newly developed technique has been proven solid against high nonlinearities in the dynamics. However, the accuracy of the KO control solution comes at a cost of few limitations and drawbacks. Due to the similarities of the KO framework to perturbation theory~\cite{servadio2021dynamics}, the robustness of the solution decreases as the relative distance between chaser and target increases. This drawback is overcome in two separate ways. First, for extremely large initial conditions, the linear dynamics solution can be used as a first approach to the desired state, switching to the high-order solution once closer. Secondly, distances can be furthermore scaled, such that the KO solves for the costates as the initial relative position is closer to the desired final state. The optimal values of the costates for the unscaled system are determined by multiplying the scaled costates by the inverse of the scaling factor. These techniques are often required to let the KO work with the scaled states, which are bounded in the [-1,1] domain or close to the origin, such that the Galerkin method best evaluates the projection of the dynamics onto the Koopman eigenfunctions. 

Lastly, the Koopman approximation of the solution of the system can only be so accurate up to certain time, before diverging from the actual true pathway. Therefore, the proposed control technique is not suitable for extremely large time constraints, like multiple weeks in the high-order Clohessy-Wiltshire equations.

%%%%%%%%%%%%%%%%%%%%%%%%%%%%%

\section{Conclusion}
A novel energy optimal inverse control technique has been presented in this paper. The state of the system is augmented to comprehend the costate variables and the resulting augmented system of ordinary differential equation is solved in the Koopman framework. Thus, the Koopman solution approximates the time evolution of the state as a linear combination well-selected basis functions: the Legendre polynomials have been selected in this work. The Koopman propagation of the state constitutes the polynomial transition map of the system. This map can be inverted to evaluate and evaluated to obtain the initial values of the optimal costates, which lead the system to the desired state at the final time. The solution of the TPBVP is therefore achieved considering the separate influence of the constraints in the Koopman solution: the time requirement is imposed by fixing the coefficients of the polynomials, while the state requirement is solved by evaluation the Legendre polynomials at the initial conditions.

The application to the Duffing oscillator offers the explanation of the new theory in a simple toy problem, and it gives a nice visualization of the solution. Afterwards, the Clohessy-Wiltshire problem has been selected to achieve rendezvous between two satellites  orbiting in close proximity. The classic linear system of equations has been presented as a particular case from the high-order Lagrangian, where the polynomial series expansion of the potential energy of the system has been truncated at the second order. Therefore, a comparison between the controlled high-order system, where the nonlinearities have been taken onto account, with the controlled linear equations, shows the improved accuracy of the dynamics. The addition of the nonlinear terms accounts for the drift of the linear periodic solution form the $y$ axis and, consequently, the optimal costate values calculated with a large truncation order are more truthful of the real behaviour of the system. The results show that the KO control technique can evaluate costates inverting a polynomial map from a Lagrangian composed by the first six polynomials of the potential. The robustness of the KO technique is tested to prove that the chaser can approach the target from any direction.

%%%%%%%%%%%%%%%%%%%%%%%%%%%%%%

\section*{Acknowledgments}
The authors want to acknowledge the support of this work by the Air Force’s Office of Scientific Research under Contract Number FA9550-18-1-0115. 

\bibliography{biblio.bib}

\end{document}